# Urban Logistics in Amsterdam: A Modal Shift from Roadways to Waterways

Nadia Pourmohammad-Zia[a], Mark van Koningsveld[a,b]

[a]Department of Hydraulic Engineering, Delft University of Technology, The Netherlands

[b]Van Oord Dredging and Marine Contractors, 3068 NH Rotterdam, The Netherlands

**Abstract–** The efficiency of urban logistics is vital for economic prosperity and quality of life in cities. However, rapid urbanization poses significant challenges, such as congestion, emissions, and strained infrastructure. This paper addresses these challenges by proposing an optimal urban logistic network that integrates urban waterways and last-mile delivery in Amsterdam. The study highlights the untapped potential of inland waterways in addressing logistical challenges in the city center. The problem is formulated as a two-echelon location routing problem with time windows, and a hybrid solution approach is developed to solve it effectively. The proposed algorithm consistently outperforms existing approaches, demonstrating its effectiveness in solving existing benchmarks and newly developed instances. Through a comprehensive case study, the advantages of implementing a waterway-based distribution chain are assessed, revealing substantial cost savings (approximately 28%) and reductions in vehicle weight (about 43%) and travel distances (roughly 80%) within the city center. The incorporation of electric vehicles further contributes to environmental sustainability. Sensitivity analysis underscores the importance of managing transshipment location establishment costs as a key strategy for cost efficiencies and reducing reliance on delivery vehicles and road traffic congestion. This study provides valuable insights and practical guidance for managers seeking to enhance operational efficiency, reduce costs, and promote sustainable transportation practices. Further analysis is warranted to fully evaluate the feasibility and potential benefits, considering infrastructural limitations and canal characteristics.

## 1. Introduction

Efficient urban logistics can be seen as a fundamental prerequisite for the economy and livability of the cities. The ever-increasing population in urban areas puts this efficiency under pressure and forces serious challenges such as congestion, emissions, noise, and safety issues on the other hand. Accordingly, seeking innovative solutions to mitigate the adverse effects of urban logistics and improve its performance is imperative. One such solution lies in the exploration of alternative transportation modes. While road transport currently dominates urban freight, inland waterways present untapped potential in many cities.

In Amsterdam, the growing strain on public spaces, along with congestion and the considerable task of maintaining bridges and quay walls, has prompted the need to reassess the city's current logistics design. The situation is particularly challenging in the historic center, where freight transport contributes to quay wall deterioration and congestion on the narrow roads alongside the canals of Amsterdam (Korff et al., 2022). Based on the topological features of the city center, located on the bank of canals, inland waterways can play a vital role in addressing its logistical needs. Currently, very few initiatives incorporate these canals for freight, indicating that there is still capacity within canals for a shift from roadways to waterways.

As highlighted by the municipality (Nepveu & Nepveu, 2020), the canals of Amsterdam present potential opportunities for a modal shift in three key areas: construction material, waste, and food supply flows. Construction material flows are typically project-based, with temporary barges serving as transshipment points at construction sites. For waste, specific transshipment infrastructure is generally unnecessary, and transfer to vessels occurs through dumping. However, for food supplies, dedicated transshipment locations are required, where goods are unloaded from vessels and loaded onto light electric vehicles for last-mile delivery. The lack of these transshipment locations poses a significant obstacle to achieving a modal shift in Amsterdam. The establishment of such points is influenced not only by considerations of subsequent routing decisions for vessels and last-mile delivery but also by the availability of space, the condition of surrounding infrastructure such as quay walls and bridges, and the characteristics of the canal classes.

Driven by the profound impact of efficient urban logistics on the overall well-being and development of cities, this study designs an optimal HoReCa logistics network for the historical center of Amsterdam. The proposed network comprises a synergistic integration of urban waterways and last-mile delivery via road transportation. To this end, the problem is formulated as a two-echelon location routing problem with time windows, where heterogeneous vessels in the first echelon and moving jacks together with light electric vehicles in the second echelon are applied. To tackle this intricate problem, a hybrid solution approach is devised, leveraging a custom-designed Adaptive Large Neighborhood Search (ALNS), local search techniques, K-means clustering, and branch and price methods. This comprehensive approach is adept at solving medium to large-sized instances of the problem, facilitating effective decision-making in urban logistics optimization for Amsterdam.

Amsterdam is not the only city facing the challenges outlined. Many cities worldwide share the need for urgent action to safeguard their historical infrastructure while enhancing their logistical networks. The delicate balance between conserving heritage and accommodating new demands is a shared concern. Historically, canals served as vital transportation routes in many city centers across the globe, leaving behind remnants of old quay walls and bridges. By learning from Amsterdam's experience and exploring innovative approaches, cities globally can strive to find sustainable solutions that enhance both their cultural legacy and efficient transportation systems.

The research on waterways for urban logistics is highly confined, and when it comes to operational and tactical decisions, this deficiency is even more highlighted. As such, a significant contribution of this research lies in developing an optimal logistics chain specifically tailored for urban waterway distribution. In this respect, we take canal classes into account, which impacts deriving the distances between different nodes for different vessel sizes and, thereby, their routing scenarios. The other novel features of this paper are raised by the incorporated modeling assumptions, which are hardly heeded even in classic urban distribution models (see Table 1). These include the decisions on establishing transshipment points (location), allowing for applying electric vehicles, synchronization, and applying moving jacks together with light vehicles. To efficiently address this complex problem, an advanced solution algorithm has been devised, combining the branch and price technique for the first tier with an Adaptive Large Neighborhood Search (ALNS) that incorporates custom-designed destroy and repair operators. The ALNS employs local search and K-means clustering techniques to enhance solution intensification in the second tier.

The subsequent sections of this paper are structured as follows: In Section 2, the relevant literature is reviewed, allowing for the identification of existing research gaps. Section 3 provides the problem description and the mathematical model. The solution methodology employed is elaborated upon in Section 4. Section 5 encompasses a detailed analysis of numerical results, the case study, sensitivity analysis, and discussion. Finally, in Section 6, the paper is concluded, emphasizing key findings and offering suggestions for potential future research avenues.

## 2. Literature Review

In this section, we briefly review the existing literature on urban logistics, where operational planning and logistics chain design are particular points of interest. Our research mainly builds on two streams: the application of inland waterways in urban logistics and two-echelon routing in urban logistics, each of which will be reviewed as follows.

### 2.1. Application of Inland Waterways in Urban Logistics

Despite its potential, the literature on waterborne freight transport in urban logistics is fairly confined, and there are limited research works in this area.

Janjevic and Ndiaye (2014), Maes et al. (2015), Miloslavskaya et al. (2019), and Wojewódzka-Król and Rolbiecki (2019) have reviewed successful practices of waterborne urban logistics worldwide. Several initiatives have been introduced in these papers; among those are:
- Beer Boat in Utrecht for deliveries to cafés, hotels, restaurants
- "Vracht door de gracht" (freight through canals) by Mokum Maritiem in Amsterdam for deliveries to local shops and waste collection
- The DHL floating service center in Amsterdam for parcel delivery
- Vert Chez Vous in Paris for parcel delivery
- Sainsbury's in London for deliveries to supermarkets
- POINT-P in Paris for construction material
- Franprix in Paris for deliveries to supermarkets
- Domestic waste transport in Lille
- Paper recycling by Barge in Paris
- Waste transport by barge in Tokyo

Kortmann et al. (2018) investigated the potential of waterborne distribution for same-day delivery to shopkeepers in Amsterdam. They developed a simulation model to analyze the performance of this distribution system and determine the appropriate fleet size. Their results show that waterborne distribution with few hubs can be a competent and sustainable delivery mode in Amsterdam, provided that further studies on its financial viability are carried out.

Gu and Wallace (2021) developed one of the few optimization models for waterborne urban logistics where the application of autonomous vessels is investigated. Their mixed-integer programming model tackles the facility location, fleet allocation, and routing decisions in the daily operations of vessels in Bergen, Norway. While their results demonstrate the potential benefits of autonomous vessels at the operational level, further investigations are required on the impact of initial investment costs.

Divieso et al. (2021) studied the feasibility of using waterways for urban logistics in Brazil. They identified features of waterway urban logistics practices worldwide and then comparatively analyzed the case for a city (Belém) in Brazil. Their analysis highlights the great potential of Belém for the application of waterways as an aid to urban logistics. Nepveu and Nepveu (2020) assessed the potential of implementing urban waterway transport in Amsterdam by exploring the success and failure factors for such a modal shift.

## 2.2. Two-Echelon Routing Problem in Urban Logistics

Two-Echelon Vehicle Routing Problem (2E-VRP) is an affluent area of academic research in urban logistics. Intermediate facilities, known as transshipment points or satellites applied for consolidation and transshipment of the items between the two echelons, are an essential part of two-echelon networks. When these points are not pre-established, and their locations need to be determined, the problem turns into a Two-Echelon Location Routing Problem (2E-LRP). Other variants of the problem, such as two-echelon inventory routing, truck-and-trailer routing, and production routing, are not the point of our interest in this paper.

Crainic et al. (2009) introduced the first 2E-VRP model in a multi-product and multi-depot setting. Zhou et al. (2018) investigated a multi-depot 2E-VRP with delivery options, allowing the customers to pick up their parcels at satellites. They developed a hybrid multi-population genetic algorithm to solve the problem. Belgin et al. (2018) studied a variant of the problem for which pick-up and deliveries are considered and applied this two-echelon distribution system in a supermarket chain in Turkey.

Li et al. (2020) investigated the application of Unmanned Aerial Vehicles (UAVs) in the second echelon, where the first echelon distribution vans are considered mobile satellites for UAVs. Enthoven et al. (2020) introduced covering locations in 2E-VRPs, from where customers can pick up their parcels. Similarly, Vincent et al. (2021) designed a two-echelon distribution system considering covering locations and occasional drivers. They showed that using crowds as occasional drivers in addition to the city freighters increases the efficiency of the distribution network.

Synchronizing the arrival and departure of the vehicles in the first and second echelons is essential in designing a seamless distribution network. Yet, this is mostly overlooked in the existing literature. Anderluh et al. (2021) provided one of the few works that took satellite synchronization into account. They considered synchronization in a multi-objective setting. In order to address the desires of citizens and municipalities, they applied a second objective function that accounts for the negative effects of transport, such as emissions. Li et al. (2021) and Jia et al. (2022) are the other researchers who took synchronization into account in their classic 2E-VRP.

The decision to establish transshipment points, leading to 2E-LRPs, increases the complexity of the already complex 2E-VRPs. Thereby, few papers have taken both locating and routing decisions into account. Among those are Zhao et al. (2018), Darvish et al. (2019), and Mirhedayatian et al. (2021), who have investigated 2E-LRP in capacitated, timely-flexible and synchronized settings, respectively. Two-echelon electric vehicle routing is the other extension of 2E-VRP, where the vehicles (mostly in the second echelon) have a limited driving range. This assumption is heeded in three ways: limiting the traveled distance or considering refilling batteries at charging or battery swap stations. Breunig et al. (2019), Jie et al. (2019), and Wu and Zhang (2021) are examples of this extension.

Developing efficient solution approaches to solve the well-established 2E-VRP and its variants is an active research direction in the area. The idea is to obtain better solutions for the existing benchmark instances and improve the gaps. Adaptive Large Neighborhood Search (ALNS) (Grangier et al., 2016), Memetic Algorithm (MA) (Bevilaqua et al., 2019), Branch and Price (B&P) (Mhamedi et al., 2022), Variable Neighborhood Search (VNS) (Akbay et al., 2022) and other heuristics such as a novel Construction Heuristic (CH) (Yu et al., 2020), Graph-Guided Heuristic (GGH) (Huang et al., 2021), Sample Average Approximation (SAA) technique (Pina-Pardo et al., 2022), and Cluster-based Artificial Immune Algorithm (C-AIA) (Liu et al., 2023) are among applied methods. Table 1 provides a general overview of the existing literature on the two-echelon routing problem in urban logistics.

Coming up with a general overview, the exploration of waterways as a viable option for urban logistics has garnered limited attention within the research community. Existing studies primarily focus on examining successful use cases rather than delving into comprehensive investigations. On the other hand, 2E-VRP in urban logistics is a relatively affluent area of academic research. However, amidst this rich tapestry of research, several promising avenues have remained undervalued or overlooked. These areas encompass the identification of optimal transshipment points, satellite synchronization for improved efficiency, integration of electric vehicles into the logistics network, consideration of multiple delivery modes, and harnessing the potential of waterways in the primary echelon of transportation. By shedding light on these neglected aspects, we can uncover new insights and opportunities to pave the way for a more sustainable and efficient urban logistics ecosystem.

## 3. Problem Description and Mathematical Model

This paper considers a multi-modal two-echelon location and routing problem with time windows that rises in urban logistics. The network embraces a combination of inland waterways and streets. The first echelon involves the flow of inland vessels from a central hub to transshipment locations in the city center and then back to the hub. The transshipment locations are the points where the vessels are unloaded, and the last-mile delivery starts by Light Electric Vehicles (LEVs) with a maximum weight of 700 kg or moving jacks. Each vessel can serve several transshipment locations, and each transshipment location can be visited by more than one vessel, implying that split delivery is admissible in the first echelon. The second echelon includes the flow of LEVs from vehicle depots to the transshipment locations and then navigating a prescribed route to serve designated demand points (HoReCa businesses) and finally returning to the depot. The demand

Table 1. General overview of the existing literature on 2E-VRP in urban logistics

| Reference | Location | Electric Vehicles | Satellite Capacity | Time Windows | Synchronization | Multiple Delivery Modes | Waterways | Solution Approach |
|---|---|---|---|---|---|---|---|---|
| Crainic et al. (2009) | | | ✓ | ✓ | | | | HD |
| Grangier et al. (2016) | | | | ✓ | | | | ALNS |
| Belgin et al. (2018) | | | ✓ | | | | | VND+LS |
| Zhao et al. (2018) | ✓ | | ✓ | | | | | CAH |
| Zhou et al. (2018) | | | ✓ | | | ✓ | | HMPG |
| Bevilaqua et al. (2019) | | | | | | | | MA |
| Breunig et al. (2019) | | ✓ | ✓ | | | | | LNS |
| Darvish et al. (2019) | ✓ | | ✓ | | | | | B&P |
| Jie et al. (2019) | | ✓ | ✓ | | | | | ALNS + B&P |
| Enthoven et al. (2020) | | | | | | ✓ | | ALNS |
| Li et al. (2020) | | | ✓ | ✓ | | | | ALNS |
| Yu et al. (2020) | | ✓ | | ✓ | | | | CH+LNS |
| Anderluh et al. (2021) | | | | | ✓ | | | LNS |
| Huang et al. (2021) | | | ✓ | | | | | GGH |
| Li et al. (2021) | | | | ✓ | ✓ | | | ALNS |
| Mirhedayatian et al. (2021) | ✓ | | ✓ | ✓ | ✓ | | | DBH |
| Vincent et al. (2021) | | | | ✓ | | ✓ | | ALNS |
| Wu and Zhang (2021) | | ✓ | | | | | | B&P |
| Jia et al. (2022) | | | ✓ | ✓ | ✓ | | | ALNS |
| Mhamedi et al. (2022) | | | | ✓ | | | | B&P |
| Akbay et al. (2022) | | | ✓ | ✓ | | | | CH |
| Pina-Pardo et al. (2022) | ✓ | | ✓ | | | | | SSA |
| Liu et al. (2023) | | | | | | | | C-AIA |
| **This Study** | ✓ | ✓ | ✓ | ✓ | ✓ | ✓ | ✓ | ALNS +LS+ B&P |

\* HD: Hierarchical Decomposition; ALNS : Adaptive Large Neighborhood Search; B&P: Branch and Price; VND: Variable Neighborhood Descent (VND); LS: Local Search; HMPG: Hybrid Multi-Population Genetic; CAH: Cooperative Approximation Heuristic; MA: Memetic Algorithm; LNS: Large Neighborhood Search; B&B: Branch and Bound; GGH: Graph-Guided Heuristic; DBH: Decomposition-Based Heuristic; CH: Construction Heuristic; SSA: sample average approximation; C-AIA: cluster-based artificial immune algorithm

points proximate enough to the transshipment locations are served by moving jacks instead of LEVs.

The problem is modeled on a directed graph $G(V, E)$, where $V$ represents the set of vertices and $E$ is the set of arcs. $V = V_1 \cup V_2$ includes the set of vertices in the first echelon ($V_1$) and second echelon ($V_2$). The set $V_1 = \{CH\} \cup TP$ involves the central hub ($CH$) and the transshipment locations ($TP$). The set $V_2 = VD \cup TP \cup HRC$ is comprised of the vehicle depots ($VD$), the transshipment locations ($TP$), and the demand points ($HRC$). $E = \{(i,j)| i,j \in V, i \neq j, (i,j) \in A \cup B\}$ where $A$ and $B$ are the sets of admissible arcs for the first and second echelon, respectively. $A = A_1 \cup A_2$ and $B = B_1 \cup B_2 \cup B_3$ where:

$A_1 = \{(i,j)| i = CH, j \in TP\}, A_2 = \{(i,j)| i \in TP, j \in TP \cup \{CH\}\}$

$B_1 = \{(i,j)| i \in VD, j \in TP\}, B_2 = \{(i,j)| i \in TP, j \in HRC\}, B_3 = \{(i,j)| i \in HRC, j \in HRC \cup VD\}$

It should be noted that the distance between any two nodes in the first echelon, is not driven only based on the shortest path method but concerning different canal classes. Thereby, the distance between two identical nodes can differ for various vessel types with different sizes. Figure (1) illustrates a typical solution on the described graph.

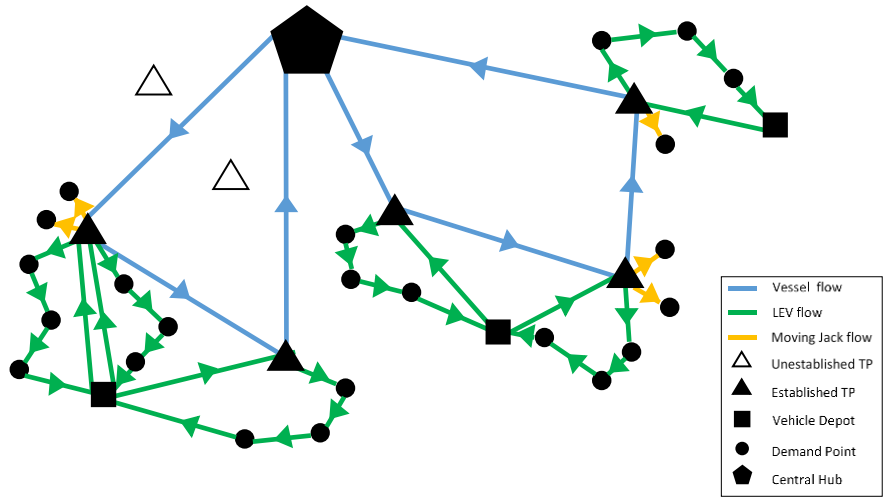

Figure 1. The graph of the problem

In order to complete the flow in the first echelon, we need to locate the transshipment points. These locations are specified from a set of potential sites for transshipment. The transshipment locations are assumed to be heterogeneous, implying that their establishment cost, capacity, and allowed laying time are different. Having the set of designated locations that are used by the vessels, LEVs, and moving jacks, the routing decisions of the vessels and LEVs are determined. The remainder of the notations which are used to formulate the model are as follows:

*Parameters*

    $D_i$     Demand of $i \in HRC$

    $S_i^{ID}$     Service time of vertex $i$ (*ID =1* first echelon, *ID=2* second echelon)

| | |
|---|---|
| $T^I_{ijk}$ | Travel time of arc (i,j) for vessel k |
| $T^{II}_{ijk}$ | Travel time of arc (i,j) for LEV k |
| $T^{III}_{ij}$ | Travel time of arc (i,j) for moving jacks |
| $AL_i$ | Allowed laying time at transshipment point i |
| $TA_i$ | Lower bound for admissible service time at vertex i |
| $TB_i$ | Upper bound for admissible service time at vertex i |
| $CAP_i$ | Capacity of transshipment point $i \in TP$ |
| $Q^I_k$ | Capacity of vessel k |
| $Q^{II}_k$ | Capacity of LEV k |
| $DL$ | Driving range limit for LEVs |
| $C_{ijk}$ | Cost of traveling arc (i,j) by vehicle k |
| $DIS_{ij}$ | Average traveling distance of arc (i,j), $i,j \in V_2$ |
| $DTr$ | Threshold distance to use moving jacks |
| $FC_i$ | Period equivalent fixed cost of establishing transshipment point $i \in TP$ |
| $l_{ij}$ | 1: if demand point $j \in HRC$ is located at a distance shorter than $DTr$ from point $i \in TP$<br>0: otherwise |
| $m_1, \ldots, m_6$ | Lower bound for the left-hand side of the respective constraints |
| $M_1, \ldots, M_4$ | Upper bound for the left-hand side of the respective constraints |

*Decision variables*

| | |
|---|---|
| $x^I_{ijk}$ | 1: if vessel $k \in K_1$ travels from $i \in V_1$ to $j \in V_1$<br>0: otherwise |
| $x^{II}_{ijk}$ | 1: if LEV $k \in K_2$ travels from $i \in V_2$ to $j \in V_2$<br>0: otherwise |
| $u_{ij}$ | 1: if the demand point $j \in HRC$ is served by a moving jack from transshipment point $i \in TP$<br>0: otherwise |
| $y_i$ | 1: if transshipment point $i \in TP$ is established<br>0: otherwise |
| $p_{ijk}$ | 1: if the items of the demand point $j \in HRC$ are delivered by vessel $k \in K_1$ to $i \in TP$<br>0: otherwise |
| $v_{k\hat{k}i}$ | 1: if LEV $\hat{k} \in K_2$ meets vessel $k \in K_1$ at $i \in TP$<br>0: otherwise |
| $st^I_{ik}$ | Time when vessel $k \in K_1$ starts to service vertex $i \in V_1$ |
| $st^{II}_{ik}$ | Time when LEV $k \in K_2$ starts to service vertex $i \in V_2$ |
| $st^{III}_{ij}$ | Time when a moving jack assigned to serve $j \in HRC$ starts to service vertex $i \in TP \cup HRC$ |
| $at^I_{ik}$ | Time when vessel $k \in K_1$ arrives at vertex i |
| $q_{ik}$ | The amount delivered by vessel $k \in K_1$ to the transshipment point $i \in TP$ |
| $pu_{ijk}$ | Auxiliary binary variable |

**Optimization Model**

$$P_1: \min Z = \sum_{k \in K_1} \sum_{i \in V_1} \sum_{j \in V_1} C^I_{ijk} x^I_{ijk} + \sum_{k \in K_2} \sum_{i \in V_2} \sum_{j \in V_2} C^{II}_{ijk} x^{II}_{ijk} + \sum_{i \in TP} FC_i\, y_i \quad (1)$$

s.t.

$$\sum_{j \in TP} x^I_{ijk} \leq 1 \qquad \forall i \in \{CH\}, k \in K_1 \quad (2)$$

$$\sum_{i \in V_1} x^I_{ivk} - \sum_{j \in V_1} x^I_{vjk} = 0 \qquad \forall v \in TP, k \in K_1 \quad (3)$$

$$\sum_{i \in V_1} x^I_{ijk} \leq y_j \qquad \forall j \in TP, k \in K_1 \quad (4)$$

$$at^I_{jk} = \sum_{i \in V_1} (st^I_{ik} + S^I_i + T^I_{ijk}) x^I_{ijk} \qquad \forall j \in V_1, k \in K_1 \quad (5)$$

$$st^I_{ik} \geq at^I_{ik} \qquad \forall i \in V_1, k \in K_1 \quad (6)$$

$$st^I_{ik} + S^I_i - at^I_{ik} \leq AL_i \qquad \forall i \in TP, k \in K_1 \quad (7)$$

$$\sum_{i \in TP} q_{ik} \leq Q^I_k \qquad \forall k \in K_1 \quad (8)$$

$$q_{jk} \leq M_1 \sum_{i \in V_1} x^I_{ijk} \qquad \forall j \in TP, k \in K_1 \quad (9)$$

$$q_{ik} = \sum_{j \in HRC} D_j\, p_{ijk} \qquad \forall i \in TP, k \in K_1 \quad (10)$$

$$\sum_{k \in K_1} \sum_{i \in TP} p_{ijk} = 1 \qquad \forall j \in HRC \quad (11)$$

$$p_{ijk} \leq \sum_{v \in V_1} x^I_{vik} \qquad \forall i \in TP, k \in K_1 \quad (12)$$

$$\sum_{k \in K_1} \sum_{j \in HRC} D_j\, p_{ijk} \leq CAP_i \qquad \forall i \in TP \quad (13)$$

$$DTr - DIS_{ij} \leq M_2 l_{ij} \qquad \forall i \in TP, j \in HRC \quad (14)$$

$$DTr - DIS_{ij} \geq m_1 (1 - l_{ij}) \qquad \forall i \in TP, j \in HRC \quad (15)$$

$$y_i + l_{ij} \leq 1 + u_{ij} \qquad \forall i \in TP, j \in HRC \quad (16)$$

$$y_i + l_{ij} \geq 2 u_{ij} \qquad \forall i \in TP, j \in HRC \quad (17)$$

$$\sum_{i \in VD} \sum_{j \in V_2} x^{II}_{ijk} \leq 1 \qquad k \in K_2 \quad (18)$$

$$\sum_{i \in V_2} x^{II}_{ijk} \leq y_j \qquad \forall j \in TP, k \in K_2 \quad (19)$$

$$\sum_{i \in V_2} x^{II}_{ivk} - \sum_{i \in V_2} x^{II}_{vjk} = 0 \qquad \forall v \in TP \cup HRC, k \in K_2 \quad (20)$$

$$\sum_{k \in K_2} \sum_{i \in V_2} x^{II}_{ijk} + \sum_{i \in V_2} u_{ij} = 1 \qquad \forall j \in HRC \quad (21)$$

$$\sum_{i \in V_2} \sum_{j \in V_2} DIS_{ij}\, x^{II}_{ijk} \leq DL \qquad \forall k \in K_2 \qquad (22)$$

$$\sum_{i \in V_2} \sum_{j \in V_2} D_j\, x^{II}_{ijk} \leq Q^{II}_k \qquad \forall k \in K_2 \qquad (23)$$

$$st^{II}_{ik} \geq (st^{II}_{ik} + S^{II}_i + T^{II}_{ijk}) x^{II}_{ijk} \qquad \forall i,j \in V_2, k \in K_2 \qquad (24)$$

$$st^{III}_{jj} = (st^{III}_{ij} + S^{II}_i + T^{III}_{ij}) u_{ij} \qquad \forall i \in TP, j \in HRC \qquad (25)$$

$$st^{II}_{i\hat{k}} - st^{I}_{ik} - S^{I}_i \geq m_2 (1 - v_{k\hat{k}i}) \qquad \forall i \in TP, k \in K_1, \hat{k} \in K_2 \qquad (26)$$

$$st^{III}_{ij} - st^{I}_{ik} - S^{I}_i \geq m_3 (1 - pu_{ijk}) \qquad \forall i \in TP, j \in HRC, k \in K_1 \qquad (27)$$

$$p_{ijk} + u_{ij} \leq 1 + pu_{ijk} \qquad \forall i \in TP, j \in HRC, k \in K_1 \qquad (28)$$

$$p_{ijk} + u_{ij} \geq 2 pu_{ijk} \qquad \forall i \in TP, j \in HRC, k \in K_1 \qquad (29)$$

$$v_{k\hat{k}j} \leq \sum_{i \in V_1} x^{I}_{ijk} \qquad \forall j \in TP, k \in K_1, \hat{k} \in K_2 \qquad (30)$$

$$\sum_{k \in K_1} v_{k\hat{k}i} = \sum_{j \in V_2} x^{II}_{ij\hat{k}} \qquad \forall i \in TP, \hat{k} \in K_2 \qquad (31)$$

$$\sum_{k \in K_1} \sum_{i \in TP} v_{k\hat{k}i} \leq 1 \qquad \forall \hat{k} \in K_2 \qquad (32)$$

$$\sum_{\hat{k} \in K_2} v_{k\hat{k}i} + u_{ij} \geq p_{ijk} \qquad \forall i \in TP, j \in HRC, k \in K_1 \qquad (33)$$

$$TA_j \sum_{i \in V_2} x^{II}_{ijk} \leq st^{II}_{jk} \leq TB_j \sum_{i \in V_2} x^{II}_{ijk} \qquad \forall j \in HRC, k \in K_2 \qquad (34)$$

$$TA_j \sum_{i \in TP} u_{ij} \leq st^{III}_{jj} \leq TB_j \sum_{i \in TP} u_{ij} \qquad \forall j \in HRC \qquad (35)$$

$$x^{I}_{ijk}, x^{II}_{ijk}, y_i, p_{ijk}, v_{k\hat{k}i}, u_{ij}, pu_{ijk} \in \{0,1\} \qquad \forall i,j \in V, k \in K \qquad (36)$$

$$st^{I}_{ik}, st^{II}_{ik}, st^{III}_i, q_{ik} \geq 0 \qquad \forall i,j \in V, k \in K \qquad (37)$$

The objective function (1) minimizes the total cost, including the travel cost of both echelons and period equivalent establishment cost of transshipment points.

Constraints (2) ensure that each vessel leaves the central hub at most once. Constraints (3) are flow constraints in the first echelon. Constraints (4) guarantee that a transshipment location can only be visited if that point is established. Consistency of the arrival time and service time of vessels is guaranteed by constraints (5) and (6). Admissible laying times are respected by constraints (7). Constraints (8) ensure that the capacity limits of the vessels are considered.

Constraints (9) express that the delivery volume of a specific vessel to a particular TP can be non-zero only if that vessel visits that TP. Quantities delivered to a TP should satisfy the demand of the HoReCa businesses allocated to this point. This is addressed by constraints (10). Each demand point is served by one TP that is guaranteed by constraints (11). Constraints (12) ensure that a TP can serve a demand point only if it is already established. Capacity limits of transshipment locations are met by constraints (13).

If a demand point is allocated to a TP located in its proximity (the distance between two is less than a pre-specified threshold), that point is served by a moving jack instead of an LEV. This is illustrated through constraints (14)-(17). Constraints (18) guarantee that each LEV can leave one of the vehicle depots and at most once. An LEV can enter a TP if that point is already established. This is addressed by constraints (19). Constraints (20) are flow constraints in the second echelon. Constraints (21) guarantee that each demand point is served either by an LEV or a moving jack. Constraints (22) respect the limited driving range of LEVs. As these vehicles are of small capacity, en-route charging is not logical for them, and battery refilling takes place at depots. Constraints (23) ensure that the capacity limits of LEVs are respected.

Consistency of t service time of the LEVs is guaranteed by constraints (24). Constraints (25) ensure the same logic for moving jacks. Constraints (26)-(33) are synchronization constraints. Expressly, constraints (26) indicate that if a vessel and an LEV are synchronized, then the service time consistency should be met. Constraints (27)-(29) express the same logic for moving jacks. Constraints (30)-(33) illustrate how the synchronization of a vessel and an LEV takes place. Time windows are represented by constraints (34) and (35). Finally, constraints (36) and (37) imply the type of variables.

Equations (5), (24), and (25) are non-linear and are linearized as follows:

$$at_{jk}^I - st_{ik}^I - S_i - T_{ijk}^I \geq m_4(1 - x_{ijk}^I) \qquad \forall i,j \in V_1, k \in K_1 \qquad (38)$$
$$at_{jk}^I - st_{ik}^I - S_i - T_{ijk}^I \leq M_3(1 - x_{ijk}^I) \qquad \forall i,j \in V_1, k \in K_1 \qquad (39)$$
$$st_{jk}^{II} - st_{ik}^{II} - S_i - T_{ijk}^{II} \geq m_5(1 - x_{ijk}^{II}) \qquad \forall i,j \in V_2, k \in K_2 \qquad (40)$$
$$st_{jj}^{III} - st_{ij}^{III} + S_i - T_{ij}^{III} \geq m_4(1 - u_{ij}) \qquad \forall i \in TP, j \in HRC \qquad (41)$$
$$st_{jj}^{III} - st_{ij}^{III} - S_i - T_{ij}^{III} \leq M_4(1 - u_{ij}) \qquad \forall i \in TP, j \in HRC \qquad (42)$$

Where constraints (38) and (39) are linearized versions of constraint (5) and constraints (40) are of constraint (24). Constraints (25) are linearized by constraints (41) and (42).

## 4. Solution Methodology

Our Two-Echelon Location Routing problem is solved by a hybrid solution algorithm that decomposes the problem into two nested sub-problems, including the first echelon and second echelon problems. We first develop an Adaptive Large Neighborhood Search (ALNS) metaheuristic to determine the location and routing decisions in the second echelon. Then, based on the provided results, we apply a Branch and Price (B&P) algorithm using the Dantzig-Wolfe decomposition principle to transform the first echelon routing model into a master problem and a subproblem. Algorithm 1. provides the pseudocode of our proposed solution methodology.

By applying the initial solution ($S_{in}$), our developed ALNS optimizes the decisions to be made for the second echelon ($S_f^2$). Then, based on the provided solution, the time windows and aggregated demand for the established points are derived. The procedure re-iterates until the termination criteria are met. The termination criteria are to reach the maximum number of iterations ($\mathbb{T}_G$) or the maximum number of iterations with no improvement.

To enhance the performance of our solution approach, we first apply preprocessing, where we remove non-admissible arcs in the second echelon concerning time windows and vehicle capacity constraints.

**Algorithm 1.** The Two-Echelon Location Routing Algorithm

**Input:** Input Parameters Data
**Output:** Best-found feasible solution ($S_f^*$)
0      $Obj(S_f^*) = \infty$
1      Generate the initial solution for the second echelon ($S_{in}^2$)
2      Based on $S_{in}^2$ generate the initial solution for the first echelon ($S_{in}^1$)
3      $S_{in} = \{S_{in}^1, S_{in}^2\}$
4      **while** termination criteria are not met
5         $S_f^{2*} \leftarrow ALNS(S_{in}^2)$    Apply ALNS to generate the best feasible solution for the second echelon
6         $S_f^{1*} \leftarrow B\&P(S_f^{2*})$    Based on $S_f^{2*}$ generate the best feasible solution for the first echelon applying B&P
7         $S_f = \{S_f^{1*}, S_f^{2*}\}$
8         **if** $Obj(S_f) \leq Obj(S_f^*)$
9             $S_f^* \leftarrow S_f$
10        **end if**
11        $S_{in} \leftarrow S_f$
12     **end while**
13     **return** $S_f^*$

### 4.1. Feasibility and Penalty Calculation

Our developed ALNS allows infeasible solutions to be a part of the search procedure and applies penalties for vehicle capacity, transshipment location capacity, and time windows violations. The generalized cost function of a solution $S$ is formulated as:

$$f_{gen}(S) = \text{obj} + \varrho_1 V_{Cap}^v(S) + \varrho_2 V_{Cap}^t(S) + \varrho_3 V_{TW}(S) + \varrho_4 V_{Dis}(S) \tag{43}$$

Where obj is the objective function (Eq. (43)) and the violations of vehicle capacity ($V_{Cap}^v(S)$), transshipment point capacity ($V_{Cap}^t(S)$), time windows ($V_{TW}(S)$), and driving range are scaled by the penalty weights $\varrho_1, \varrho_2, \varrho_3$ and $\varrho_4$, respectively. Every $\psi_P$ iterations, referred to as the penalty update period, the penalty weights are dynamically adjusted. If a constraint has been violated in at least $\mathcal{L}^{\psi_P}$ out of $\psi_P$ iterations, its respective penalty weight is multiplied by $\omega_i$ and it is divided by the same value if the limitation is met in at least $\mathcal{L}^{\psi_P}$ iterations. To control the magnitude of $f_{gen}$, the penalty factors are restricted to fall within the minimum and maximum values.

Efficient calculation of the changes in cost the function plays a crucial role in the performance of our solution algorithm. The changes in capacity and driving range violations are trivially obtained in $\mathcal{O}(1)$. In order to calculate the time windows violations, we incorporate the approach proposed by Nagata et al. (2010) and Schneider et al. (2013). Based on this approach, the violations in a node do not propagate to subsequent nodes of a sequence. As service time at different vertices is independent of the route sequence, the changes in time windows violations can be calculated in $\mathcal{O}(1)$.

### 4.2. Initial Solution

A route in the second echelon starts from a vehicle depot, heads to a transshipment location, then navigates through several demand points, and finally returns to the initial vehicle depot. A solution for the second echelon is comprised of several routes, and each route is represented by a series of points $< VD_n, TP_m, HRC_1, \ldots, HRC_i, \ldots, HRC_k, VD_n >$, where $VD_n$ is the starting and ending depot for the route, $TP_m$ is the selected and established transshipment point, and $HRC_1, \ldots, HRC_i, \ldots, HRC_k$ are the covered demand points within this route.

Barreto et al. (2007), Wang et al. (2013), and Akpunar and Akpinar (2021) are among the research works showing that grouping approaches have great potential in providing high-quality solutions for capacitated location routing problems. There exists several clustering algorithms in the literature such as DBSCAN, Gaussian Mixture, and K-means. Despite its straightforward approach, K-means has shown a good performance in clustering spatial data (Akpunar and Akpinar, 2021). In this respect, we apply the K-means clustering algorithm to specify the number and location of established transshipment points in our initial solution. K-means partitions our set of demand points into K initial clusters, within each of which a transshipment point is established to serve the demand points. Interested readers are referred to Likes et al. (2003) for comprehensive information on the structure and performance of the K-means algorithm.

In our algorithm, selecting the number of clusters is a crucial step. We first take K as the rough estimation of the lower bound on the number of required transshipment points by dividing the total demand of all demand points by the average capacity of transshipment locations. Cluster centers are chosen and updated from the set of potential candidates for locating transshipment points. Once demand points are assigned to these clusters, the required items' volume at each TP (cluster center) is obtained. If the capacity of at least one point is violated by over 25%, the number of clusters is increased by one, and K-means is re-iterated. This procedure continues until no more than 25% violation is observed. Once the established transshipment points are finalized, the closest depot to each point is selected as the starting and ending depot of the vehicles serving that transshipment point. Accordingly, the first two points for each route representation can be specified.

By having the established transshipment points, the demand points located within a distance of $DTr$, which will be served by moving jacks, are specified and removed from the set of uncovered demand points. In order to complete the remainder of each route, we apply the Semi-Parallel Construction (SPC) heuristic proposed by Paraskevopoulos et al. (2008). Our clustering scheme can enhance the performance of the applied semi-parallel construction heuristic. SPC is an iterative approach, and at each iteration, we first select a cluster randomly and thereby specify the first two points in the route. The set of uncovered demand points considered for insertion in this cluster is decreased from the global set to the uncovered points within this cluster and outside the cluster but located within a maximum distance from the transshipment point. Accordingly, a unique point can potentially be a part of constructed routes of more than one cluster, implying that our choice of the initially selected cluster may affect the constructed routes. In this respect, we specify the neighboring clusters starting with the initial cluster. These are the clusters whose assigned demand points (by K-means) are among the set of uncovered points of the initial randomly chosen cluster. These clusters are added to the list of our candidate clusters, and for each one, the set of within-cluster uncovered customers is specified. The procedure continues for all clusters in the candidate list until no further neighboring cluster exists.

At each iteration, an LEV for each cluster in the candidate list is taken to serve a part of uncovered demand points. Then, for each LEV, a route covering a part of uncovered demand points is constructed based on a greedy approach while taking capacity constraints into account. For diversification, we take a Restricted Candidate List (RCL) of $n$ demand points with the highest evaluation score and apply roulette wheel selection to select and insert a demand point. After constructing the routes, the best cluster-route is selected, and its covered customers are removed from the set of global uncovered demand points. The procedure re-iterates until no further uncovered demand points are left. The best vehicle-route in each iteration has the lowest Average Cost per Unit Transferred (ACUT), which divides the associated costs by the sum of the demand of visited demand points.

### 4.3. Adaptive Large Neighborhood Search

Adaptive Large Neighborhood Search (ALNS) is introduced by Ropke and Pisinger (2006a) and incorporates destroy and repair operators based on their success rates in previous iterations. We have taken ALNS as the core of our solution algorithm for the second echelon and adapted it to fit into our problem by allowing infeasible solutions, incorporating problem-specific destroy and repair operators, and applying local search for intensification. In order to enhance the speed of our algorithm, the potential demand points for insertion in repair and local search operators are selected from the neighboring clusters. That is, a demand point that is considerably far from a transshipment location cannot be served by that TP. The overview of our proposed ALNS is illustrated through Algorithm 2.

The termination criteria are to reach the maximum number of iterations ($\mathbb{T}_m$) or maximum number of iterations with no improvement ($\mathbb{T}_{n-m}$). We further apply a restarting mechanism (Hiremann et al., 2016) to escape possible local optima. More specifically, after certain iterations ($\psi_R$) with no improvement, the current solution is replaced by the best-obtained solution ($S^{2*}$). We will take a more in-depth look at different components of the proposed ALNS in the following subsections.

#### *Destroy and Repair Operators*

We have two types of destroy operators. The first, with a large impact, changes a part of the problem's configuration by removing transshipment points, and the second, with a small impact, affects only a part of the constructed routes. A destroy operator removes at least $\mathcal{G}$ nodes from the current solution, where $\mathcal{G}$ is a percentage of all nodes. This percentage is randomly drawn from a given interval $[\Lambda_{min}, \Lambda_{max}]$. The following destroy operators are used in our ALNS algorithm, where the first four are operators with small impact and the rest with large.

**Random removal** removes $\mathcal{G}$ arbitrary nodes from the current solution, where all nodes have equal removal chances.

**Algorithm 2.** ALNS

**Input:** Initial Solution ($S_{in}^2$)
**Output:** Best obtained feasible solution for the second echelon ($S_f^{2*}$)

1  $S^2 \leftarrow S_{in}^2$ and $S^{2*} \leftarrow S_{in}^2$
2  **if** $S_{in}^2$ is feasible
3  $\quad S_f^{2*} \leftarrow S_{in}^2$
4  **end if**
5  $i, j \leftarrow 0$
6  **while** termination criteria are not met
7  $\quad S^{2\prime} \leftarrow$ Destroy&Repair $(S^2)$
8  $\quad S^{2\prime} \leftarrow$ LocalSearch $(S^{2\prime})$
9  $\quad$ **if** $S^{2\prime}$ is accepted
10 $\quad\quad S^2 \leftarrow S^{2\prime}$
11 $\quad\quad$ **if** $f_{gen}(S^{2\prime}) \leq f_{gen}(S^{2*})$
12 $\quad\quad\quad S^{2*} \leftarrow S^{2\prime}$
13 $\quad\quad$ **end if**
14 $\quad\quad$ **if** $S^{2\prime}$ is feasible **and** $f_{gen}(S^{2\prime}) \leq f_{gen}(S_f^{2*})$
15 $\quad\quad\quad S_f^{2*} \leftarrow S^{2\prime}$
16 $\quad\quad\quad j \leftarrow i$
17 $\quad\quad$ **end if**
18 $\quad$ **end if**
19 $\quad$ UpdateScore($S^{2\prime}$)
20 $\quad$ **if** $0 \equiv (i+1)\ mod(\psi_P)$
21 $\quad\quad$ UpdatePenalty($S^{2\prime}$)
22 $\quad$ **end if**
23 $\quad$ **if** $0 \equiv (i - j + 1)\ mod(\psi_R)$
24 $\quad\quad S^2 \leftarrow S^{2*}$
25 $\quad$ **end if**
26 $\quad$ **if** $0 \equiv (i+1)\ mod(\psi_L)$
27 $\quad\quad$ AdaptSelectionScores()
28 $\quad$ **end if**
29 $\quad i \leftarrow i + 1$
30 **end while**
31 **return** $S_f^{2*}$

**Worst removal** was introduced by Ropke and Pisinger (2006b), where the idea is to remove the most expensive parts of the solution. The associated cost of each node is determined as the difference between the costs of the solution, with and without that node. Then, we take a Restricted Candidate List (RCL) of $n$ nodes with the highest costs and apply roulette wheel selection to select and remove a node. The procedure continues till $\mathcal{G}$ nodes are removed.

**Shaw removal** was introduced by Shaw (1997) and removes the demand points concerning their relatedness. The relatedness between any two pairs of customers is defined based on their distance, demand, and earliest service start time and is formulated as:

$$Re(i,j) = \varpi_1 \frac{Dist_{ij}}{\max_{i,j \in \vartheta_D}(Dist_{ij})} + \varpi_2 \frac{|D_i - D_j|}{\max_{i \in \vartheta_D}(D_i) - \min_{i \in \vartheta_D}(D_i)} + \varpi_3 \frac{|TA_i - TA_j|}{\max_{i \in \vartheta_D}(TA_i) - \min_{i \in \vartheta_D}(TA_i)} \quad (44)$$

Where $\varpi_1, \varpi_2$, and $\varpi_3$ are normalizing weights. The operator removes the first demand point randomly. Then, a point is randomly selected from the list of removed demand points, and the

Shaw relatedness measure is derived for any pairs containing that customer and non-removed customers. We take a Restricted Candidate List (RCL) of $n$ demand points with the largest relatedness values and apply roulette wheel selection to select and remove a demand point. The procedure continues till $\mathcal{G}$ nodes are removed.

**Route removal** has two variants: Random Route removal and Inefficient Route removal. The first selects random routes to be removed until at least $\mathcal{G}$ demand points are removed from the current solution. The second one takes a Restricted Candidate List (RCL) of $n$ routes with the largest ACUT and applies a roulette wheel selection to select and remove a route. The procedure continues until at least $\mathcal{G}$ nodes are removed.

**Transshipment point removal** is introduced by Hemmelmayer et al. (2012), where a transshipment point is chosen randomly from the list of open ones and gets closed. Therefore, all routes originating from this TP are removed, adding their covered demand points together with points served by moving jacks from that TP to the customer pool. Furthermore, a TP is randomly chosen, and if it is not already established is opened. This prevents situations in which all transshipment points would be closed. By opening the new TP, demand points to be served by moving jacks from that transshipment location are specified and removed from their current routes or customer pool.

**Transshipment point opening** is introduced by Hemmelmayer et al. (2012), where we choose a transshipment point randomly among unestablished ones and open it. Then, all demand points that can be served by moving jacks from this TP are removed from their current routes. If the number of these removed points is smaller than $\mathcal{G}$, we continue by removing the remainder of demand points from the set of closest ones to this TP.

**Transshipment point swap** is introduced by Hemmelmayer et al. (2012), where the TP removal operator is applied first. Then, we use a Restricted Candidate List (RCL) of $n$ unestablished clusters with the shortest distance from the removed TP and apply a roulette wheel selection to select and establish a new transshipment location. Similarly, By opening the new TP, demand points to be served by moving jacks from that transshipment location are specified and removed from their current routes or customer pool.

Our ALNS applies the following repair operators.

**Greedy insertion** investigates the increase in cost associated with adding each unassigned demand point to each position of a route. Then, the route and position with the lowest cost increase are selected.

**Regret insertion** was introduced by Ropke and Pisinger (2006b), where the idea is to first insert a demand point in the best route and position for which a later insertion causes the highest additional costs. The operator specifies the best position and insertion cost of all current partial routes for each demand point. A k-regret value of a customer, which projects the cost difference between the cheapest route and k-1 next cheapest routes, is calculated. Then, the demand point with the largest k-regret value is selected, and the procedure re-iterates until all demand points are inserted.

**SPC-based insertion** modifies the SPC heuristic of the construction phase by inserting demand points to existing and new routes, provided that capacity constraints are met.

*Local Search*

To further improve the results of the destroy and repair operators, a Local Search (LS) is applied that uses a composite neighborhood, including the well-known 2-opt (Potvin & Rousseau, 1995), $2^*$-opt (Potvin & Rousseau, 1995), Reinsertion (Savelsbergh, 1992), and Swap$(n-1)$, where $n = 1, \ldots, 4$ (Savelsbergh, 1992) moves.

The list of the applied neighborhoods is randomly ordered, and the algorithm starts by searching the first one until no further cost reductions can be achieved. Then, the next neighborhoods are taken and searched one after another. Once the last neighborhood of the list is searched, the procedure starts again and iteratively continues until a local minimum is obtained. To speed up the procedure, in each iteration of our LS, the first $J$ generated moves are taken, and the best-improving move is selected. If no improving move exists among these moves, the procedure continues until an improving solution is found or the neighborhood is searched completely.

*Acceptance Criteria*

We apply Simulated Annealing (SA)-based acceptance criteria in our developed ALNS: If a solution is improving, always accept it. Decide about non-improving (deteriorating) solutions based on a probability that depends on the amount of solution deterioration and the temperature following a cooling scheme.

We store two best solutions during iterations of ALNS: the best feasible solution found so far ($S_f^{2*}$) and the best non-necessarily feasible solution ($S^{2*}$) with its current penalty values. So, we need to adapt the value of the best solution after each change in penalty weights. This can lead to a $S^{2*}$ that is higher than the best feasible solution, where it should be replaced by $S_f^{2*}$.

*Adaptive Mechanism*

The destroy and repair operators are selected using a roulette wheel mechanism. The selection probability of each operator relies on its historical performance, which is projected by a weight. These weights are equal at the beginning of the algorithm and are updated after each $\psi_L$ iterations, which is referred to as the adaption period. The weight of an operator $i$ in adaption period $t$ is updated as:

$$\mathcal{W}_{it} = \theta \frac{SN_{it}}{NN_{it}} + (1-\theta)\mathcal{W}_{i\ t-1} \tag{45}$$

Where $\theta$ is the smoothing factor, $SN_i$ shows the success score of the operator in the current update period and $NN_{it}$ reflects the number of times the operator has been applied in update period $t$. In each update period, $SN_{it}$ is initially set to zero and is updated by a scoring scheme using $SN_{it} = SN_{it} + \sigma_i, i = 1,2,3$.

- $\sigma_1$ is applied when the generated solution is the best solution found so far.
- $\sigma_2$ is applied when the generated solution has not been accepted before and improves the current solution.
- $\sigma_3$ is applied when the generated solution has not been accepted before, is non-improving but accepted.

## 4.4. Branch and Price

Once the solution of the second echelon is obtained, the problem at the first echelon is seen as a split delivery vehicle routing problem with time windows, where the established transshipment points are seen as the demand points. Next, we need to specify the demand and the time windows at TPs. Since split deliveries are admissible, we cannot treat the set of all allocated demand points of a TP as a unit point. This is because all points visited by a single LEV should be served by a unique vessel due to synchronization constraints. Accordingly, we will have the accumulation of demands served by each LEV as one unique demand, located at its initial TP and with time windows respecting the time windows of all those demand points. Then, we need to create copies of TPs with demands equal to the demand of points served by a moving jack or accumulated demand of each LEV.

This potentially can lead to many serving points with the same location and thereby high degeneracy of the problem. In order to mitigate this, demand points can be merged under certain conditions. Since split delivery was to resolve the problem of limited vessel capacity and considerably different time windows, the following model can be applied to merge the points in an efficient way for each TP:

*Parameters*
- $CAP$     A percentage of the smallest vessel's capacity (e.g., 25 %)
- $TR$     The threshold for the difference in time windows

*Decision variables*
- $\mu_i$     1: if group $i \in I$ is formed; 0: otherwise
- $\lambda_{ij}$     1: if point $j \in J$ is merged into group $i \in I$; 0: otherwise

$$P_2: \min Z = \sum_{i \in I} \mu_i \tag{46}$$

s.t.

$$\sum_{j \in J} D_j \lambda_{ij} \leq CAP \cdot \mu_i \quad \forall i \in I \tag{47}$$

$$\sum_{i \in I} \lambda_{ij} = 1 \quad \forall j \in J \tag{48}$$

$$|TA_j - TA_{j'}| \lambda_{ij} \lambda_{ij'} \leq TR \,^1 \quad \forall i \in I, j, j' \in J \tag{49}$$

$$|TB_j - TB_{j'}| \lambda_{ij} \lambda_{ij'} \leq TR \quad \forall i \in I, j, j' \in J \tag{50}$$

---

[1] Constraints (49) and (50) are non-linear and can be linearized as:

$(TA_j - TA_{j'})\eta_{ijj'} \leq TR$     $\forall i \in I, j, j' \in J$
$(TA_{j'} - TA_j)\eta_{ijj'} \leq TR$     $\forall i \in I, j, j' \in J$
$(TB_j - TB_{j'})\eta_{ijj'} \leq TR$     $\forall i \in I, j, j' \in J$
$(TB_{j'} - TB_j)\eta_{ijj'} \leq TR$     $\forall i \in I, j, j' \in J$
$\lambda_{ij} + \lambda_{ij'} \leq 1 + \eta_{ijj'}$     $\forall i \in I, j, j' \in J$
$\lambda_{ij} + \lambda_{ij'} \geq 2\eta_{ijj'}$     $\forall i \in I, j, j' \in J$

$$\mu_i, \lambda_{ij} \in \{0,1\} \qquad \forall i \in I, j \in J \qquad (51)$$

The problem is a variant of the bin packing problem for which a strong valid inequality exists as follows (Correia et al., 2008).

$$\lambda_{ij} \leq \mu_i \qquad \forall i \in I, j \in J \qquad (52)$$

$$\sum_{i \in I} \mu_i \geq |I| - r + 1 \qquad \forall i \in I, j, j' \in J \qquad (53)$$

Where $|I|$ is the size of potential groups, and r is a value satisfying the following inequality:

$$(|I| - r - 1)CAP < \sum_{j \in J} D_j \leq (|I| - r)CAP \qquad \forall i \in I, j, j' \in J \qquad (54)$$

Since the capacity of transshipment points is limited, the size of $J$ for each TP is rather confined. Therefore, applying this valid inequality, the problem can be solved to optimality in a reasonable time.

In this way, the problem at the first echelon is transformed into a classic VRP with time windows, to solve which there exists extensive research applying branch and price. The approach works based on Dantzig Wolfe Decomposition (DWD), where the main optimization problem is decomposed into a master and several sub-problems (see Desaulniers et al. (2006) for details). It exploits the fact that in a classic VRP with time windows, the constraint associated with serving the demand point with one of the existing vehicles is the only constraint linking the vehicles together. By neglecting this constraint, the problem can be decomposed into sub-problems (each for one vehicle) that take the form of the shortest path problem with resource constraint (time windows and vehicle capacities).

## 5. Numerical Results

In this section, we provide the results of the conducted numerical experiments on the developed solution approach for solving our proposed variant of the two-echelon location routing problem with time windows and synchronization. As illustrated in Table 1, the problem is new, for which no benchmark instances exist. Accordingly, the experiments are carried out on our newly-generated benchmark instances. Furthermore, we assess the performance of our developed solution approach on the available benchmark instances for 2E-VRP. Finally, a case study is presented to illustrate the results of the problem in a practical setting for the city of Amsterdam, followed by sensitivity analysis on the input parameters and discussion.

The experiments are conducted on a computer with Intel® Core i7-8650U CPU 1.9 GHz, 2.11 GHz, and 32 GB memory available. Our developed solution approach was implemented and run on Python 3.6 and applied IBM ILOG CPLEX Optimization Studio 12.7.

### 5.1. Parameter Tuning and Generation of Benchmark Instances

For tuning parameters, we first searched the literature for the existing values and applied values in different reasonable ranges for the newly introduced (non-existing) parameters. To achieve a setting with a good performance, we then investigated the impact of modifying the values of these parameters on a number of problem instances, each time altering one and remaining others unchanged. The selected setting, that provided the best results found and is applied in our numerical experiments, is represented in Table 2.

Table 2. Applied parameter setting of the developed solution algorithm

| General | | ALNS | |
|---|---|---|---|
| $\mathbb{T}_G$ | 50 | $\mathbb{T}_m$ | 2000 |
| $(\varrho_1, \varrho_2, \varrho_3, \varrho_4)$ | (10,10,10,10) | $\mathbb{T}_{n-m}$ | 250 |
| $(\varrho_1^{min}, \varrho_2^{min}, \varrho_3^{min}, \varrho_4^{min})$ | (0.1,0.1,0.1,0.1) | $\psi_R$ | 200 |
| $(\varrho_1^{max}, \varrho_2^{max}, \varrho_3^{max}, \varrho_4^{max})$ | (5000,5000,5000,5000) | $\psi_L$ | 50 |
| $\psi_P$ | 10 | $[\Lambda_{min}, \Lambda_{max}]$ | [0.05,0.15] |
| $\mathcal{L}^{\psi_P}$ | $0.25 \times \psi_P$ | $(\varpi_1, \varpi_2, \varpi_3)$ | (6,4,5) |
| $(\omega_1, \omega_2, \omega_3)$ | (1.2,1.2,1.2) | $\mathcal{J}$ | 50 |
| $n$ | 5 | $\theta$ | 0.6 |

The benchmark instances involve ten instance sets classified based on their number of customers ranging from 5 to 200 and labeled as small (5-25 customers), medium (50 and 75 customers), and large (100-200 customers) size sets. The customer locations are randomly selected within a target zone, and the zone is considered to occupy an area of 0.5 km² for instances of 5, 10, 15, 20, and 25 customers, 1 km² for 50 customers, 5 km² for 75 customers, 10 km² for 100 customers, 15 km² for 150 customers, and 25 km² for 200 customers. The number of LEV depots in small instances is either 1 or 2, while for medium and large size instances, this value is 3 or 4. The number of transshipment points ranges from 2 to 4 in small size and 5 to 10 in medium and large size instances. Then, an instance with ID *SI-Dm-Cn-To* represents a small size instance with *m* depots, *n* customers, and *o* transshipment points.

Distances are transformed into travel time by considering speeds of 30 km/h for LEVs and 5 to 15 km/h for vessels traveling from the central hub to TPs. The unit energy consumption cost is 0.27 €/km for LEVs and ranges between 1.8 and 2.5 €/km for vessels. The period equivalent fixed cost of establishing transshipment points ranges from 100 to 175 €, based on the space of that TP.

## 5.2. Experiment on Existing Benchmark Instances

In order to assess the performance of our proposed ALNS+B&P, we conduct experiments on existing benchmark instances for 2E-VRP, including two instance sets (sets 2 and 3) from Perboli et al. (2011) and one larger instance set (set 5) from Hemmelmayer et al. (2012). In this regard, assumptions associated with locating TPs, electric vehicles, time windows, synchronization, and multiple delivery modes are relaxed. Therefore, all related procedures are eliminated from the solution approach, and our ALNS+B&P is decreased to solve a classic 2E-VRP.

These three sets include 21, 18, and 18 instances, where for sets 2 and 3, the number of customers ranges from 21 to 50, with 2 or 4 transshipment locations, and for set 4, the number of customers is either 100 or 200, with 5 or 10 transshipment locations. In Table 3, the performance of our proposed ALNS on these benchmark instances is compared with those of Hemmelmayer et al. (2012), Breunig et al. (2016), Enthoven et al. (2020), and Vincent et al. (2021).

The second column of the table provides the Best-Known Solution (BKS) reported in Breunig et al. (2016). The next columns present the gap of the four papers' results to this BKS, and finally, we have the results of our proposed ALNS+B&P. The average run time and average gaps are reported in the last two rows of the table. It should be noted that these algorithms were run on different platforms, and our proposed solution approach is the only one among these four that incorporates an exact algorithm (B&P) in its structure. Accordingly, a precise comparison of time is not possible.

Table 3. The results of experiments on existing 2E-VRP benchmark instances

| Instance | BKS | $\Delta_{BKS}$ | | | | |
|---|---|---|---|---|---|---|
| | | Hemmelmayr et al. (2012) | Breunig et al. (2016) | Enthoven et al. (2020) | Vincent et al. (2021) | Our proposed approach |
| **Set 2** | | | | | | |
| E-n22-k4-s6-17 | 417.07 | 0.00 % | 0.00 % | 0.00 % | 0.00 % | 0.00 % |
| E-n22-k4-s8-14 | 384.96 | 0.00 % | 0.00 % | 0.00 % | 0.00 % | 0.00 % |
| E-n22-k4-s9-19 | 470.6 | 0.00 % | 0.00 % | 0.00 % | 0.00 % | 0.00 % |
| E-n22-k4-s10-14 | 371.5 | 0.00 % | 0.00 % | 0.00 % | 0.00 % | 0.00 % |
| E-n22-k4-s11-12 | 427.22 | 0.00 % | 0.00 % | 0.00 % | 0.00 % | 0.00 % |
| E-n22-k4-s12-16 | 392.78 | 0.00 % | 0.00 % | 0.00 % | 0.00 % | 0.00 % |
| E-n33-k4-s1-9 | 730.16 | 0.00 % | 0.00 % | 0.00 % | 0.00 % | 0.00 % |
| E-n33-k4-s2-13 | 714.63 | 0.00 % | 0.00 % | 0.00 % | 0.00 % | 0.00 % |
| E-n33-k4-s3-17 | 707.48 | 0.00 % | 0.00 % | 0.00 % | 0.00 % | 0.00 % |
| E-n33-k4-s4-5 | 778.74 | 0.00 % | 0.00 % | 0.00 % | 0.05 % | 0.05 % |
| E-n33-k4-s7-25 | 756.85 | 0.00 % | 0.00 % | 0.00 % | 0.00 % | 0.00 % |
| E-n33-k4-s14-22 | 779.05 | 0.00 % | 0.00 % | 0.00 % | 0.00 % | 0.00 % |
| E-n51-k5-s2-17 | 597.49 | 0.00 % | 0.00 % | 0.00 % | 0.00 % | 0.00 % |
| E-n51-k5-s4-46 | 530.76 | 0.00 % | 0.00 % | 0.00 % | 0.00 % | 0.00 % |
| E-n51-k5-s6-12 | 554.81 | 0.00 % | 0.00 % | 0.00 % | 0.04 % | 0.00 % |
| E-n51-k5-s11-19 | 581.64 | 0.00 % | 0.00 % | 0.00 % | 0.00 % | 0.00 % |
| E-n51-k5-s27-47 | 538.22 | 0.00 % | 0.00 % | 0.00 % | 0.00 % | 0.00 % |
| E-n51-k5-s32-37 | 552.28 | 0.00 % | 0.00 % | 0.00 % | 0.00 % | 0.00 % |
| E-n51-k5-s2-4-17-46 | 530.76 | 0.00 % | 0.00 % | 0.00 % | 0.00 % | 0.00 % |
| E-n51-k5-s6-12-32-37 | 531.92 | 0.00 % | 0.00 % | 0.00 % | 0.00 % | 0.00 % |
| E-n51-k5-s11-19-27-47 | 527.63 | 0.00 % | 0.00 % | 0.00 % | 0.00 % | 0.00 % |
| **Set 3** | | | | | | |
| E-n22-k4-s13-14 | 526.15 | 0.00 % | 0.00 % | 0.00 % | 0.00 % | 0.00 % |
| E-n22-k4-s13-16 | 521.09 | 0.00 % | 0.00 % | 0.00 % | 0.00 % | 0.00 % |
| E-n22-k4-s13-17 | 496.38 | 0.00 % | 0.00 % | 0.00 % | 0.00 % | 0.00 % |
| E-n22-k4-s14-19 | 498.8 | 0.00 % | 0.00 % | 0.00 % | 0.00 % | 0.00 % |
| E-n22-k4-s17-19 | 512.8 | 0.00 % | 0.00 % | 0.00 % | 0.00 % | 0.00 % |
| E-n22-k4-s19-21 | 520.42 | 0.00 % | 0.00 % | 0.00 % | 0.00 % | 0.00 % |
| E-n33-k4-s16-22 | 672.17 | 0.00 % | 0.00 % | 0.00 % | 0.00 % | 0.00 % |
| E-n33-k4-s16-24 | 666.02 | 0.00 % | 0.00 % | 0.00 % | 0.00 % | 0.00 % |
| E-n33-k4-s19-26 | 680.36 | 0.00 % | 0.00 % | 0.00 % | 0.00 % | 0.00 % |
| E-n33-k4-s22-26 | 680.36 | 0.00 % | 0.00 % | 0.00 % | 0.00 % | 0.00 % |
| E-n33-k4-s24-28 | 670.43 | 0.00 % | 0.00 % | 0.00 % | 0.00 % | 0.00 % |
| E-n33-k4-s25-28 | 650.58 | 0.00 % | 0.00 % | 0.00 % | 0.00 % | 0.00 % |
| E-n51-k5-s12-18 | 690.59 | 0.00 % | 0.00 % | 0.00 % | 0.66 % | 0.00 % |
| E-n51-k5-s12-41 | 683.05 | 0.00 % | 0.00 % | 0.00 % | 1.70 % | 0.00 % |
| E-n51-k5-s12-43 | 710.41 | 0.00 % | 0.00 % | 0.00 % | 0.00 % | 0.00 % |
| E-n51-k5-s39-41 | 728.54 | 0.00 % | 0.00 % | 0.00 % | 0.00 % | 0.00 % |
| E-n51-k5-s40-41 | 723.75 | 0.00 % | 0.00 % | 0.00 % | 0.38 % | 0.00 % |
| E-n51-k5-s40-43 | 752.15 | 0.00 % | 0.00 % | 0.00 % | 0.26 % | 0.00 % |
| **Set 5** | | | | | | |
| 100-5-1 | 1564.46 | 0.06 % | 0.00 % | 2.63 % | 0.31 % | 0.00 % |
| 100-5-1b | 1108.62 | 0.25 % | 0.00 % | 1.18 % | 1.88 % | **-0.41 %** |
| 100-5-2 | 1016.32 | 0.00 % | 0.00 % | 0.57 % | 0.50 % | 0.00 % |
| 100-5-2b | 782.25 | 0.00 % | 0.00 % | 0.00 % | 0.06 % | 0.00 % |
| 100-5-3 | 1045.29 | 0.00 % | 0.00 % | 0.07 % | 0.05 % | 0.00 % |
| 100-5-3b | 828.54 | 0.05 % | 0.00 % | 0.05 % | 0.13 % | 0.00 % |
| 100-10-1 | 1124.93 | 0.47 % | 0.05 % | 0.01 % | 0.02 % | 0.00 % |

| Instance | BKS | $\Delta_{BKS}$ | | | | |
|---|---|---|---|---|---|---|
| | | Hemmelmayr et al. (2012) | Breunig et al. (2016) | Enthoven et al. (2020) | Vincent et al. (2021) | Our proposed approach |
| 100-10-1b | 916.25 | 0.03 % | 0.00 % | 0.55 % | 0.90 % | **-0.49** % |
| 100-10-2 | 990.58 | 0.00 % | 2.18 % | 2.03 % | 2.57 % | 0.00 % |
| 100-10-2b | 768.61 | 0.00 % | 1.65 % | 0.86 % | 3.14 % | 0.00 % |
| 100-10-3 | 1043.25 | 0.00 % | 0.62 % | 0.94 % | 0.63 % | 0.00 % |
| 100-10-3b | 850.92 | 0.00 % | 0.47 % | 0.98 % | 1.05 % | **-0.14** % |
| 200-10-1 | 1556.79 | 1.11 % | 1.51 % | 0.91 % | 0.13 % | **-1.18** % |
| 200-10-1b | 1187.62 | 1.19 % | 0.33 % | 0.63 % | 0.09 % | **-0.99** % |
| 200-10-2 | 1365.74 | 0.66 % | 0.05 % | 1.71 % | 1.92 % | 0.00 % |
| 200-10-2b | 1002.85 | 0.09 % | 0.56 % | -0.02 % | 0.42 % | 0.00 % |
| 200-10-3 | 1787.73 | 0.00 % | 0.56 % | 2.79 % | 1.47 % | **-0.45** % |
| 200-10-3b | 1197.9 | 0.24 % | 0.36 % | 1.84 % | 0.95 % | 0.00 % |
| Avg. Gap | | 0.077 % | 0.155 % | 0.328 % | 0.354 % | -0.067 % |
| Avg. Time (min) | | 4.31 | 5.67 | 5.71 | 4.76 | 6.48 |

As the table represents, the average gap of our proposed solution approach is 0.14% lower than Hemmelmayr et al. (2012), 0.22% lower than Breunig et al. (2016), 0.40% lower than Enthoven et al. (2020), and 0.42% lower than Vincent et al. (2021). It should be noted that for sets 2 and 3, our results are almost identical to BKS. This is while for set 5, which embraces larger sizes, our proposed methodology improves BKS in six cases, with an average gap of -0.20 %. Considering these three sets of instances, our developed method provides better results than Hemmelmayr et al. (2012), Breunig et al. (2016), Enthoven et al. (2020), and Vincent et al. (2021) in 12, 13, 16, and 24 cases, respectively.

### 5.3. Experiment on Newly Generated Benchmark Instances

This section provides experimental results on our newly generated benchmark instances. We first use the instance sets with 5-25 customers, 1 or 2 LEV depots, and 2 to 4 transshipment points to analyze the performance of our proposed solution algorithm on small-size instances. To this end, the results of the proposed solution algorithm are compared with the optimal (global or local) solutions provided by CPLEX.

Table 4 provides the results of our experiments on small-size instances. For CPLEX, the objective function value ($Z_1$) and run time (*t*) in seconds are reported. The computing time of CPLEX is limited by 2 hours (7200 seconds). So, optimality is not guaranteed for instances that have reached this upper bound. In our proposed solution algorithm, $Z_1$ is associated with the best-found solution in ten runs of the algorithm. $\Delta_{CPLEX}$ represents the gap of the obtained objective function value to the one provided by CPLEX.

Table 4. The results of experiments on newly generated small-size benchmark instances

| Instance | CPLEX | | Proposed Solution Algorithm | | |
|---|---|---|---|---|---|
| | $Z_1$ | *t* (s) | $Z_1$ | $\Delta_{CPLEX}$ | *t* (s) |
| SI-D1-C5-T2 | 119.5287 | 1.06 | 119.5287 | 0.000 % | 5.28 |
| SI-D2-C5-T2 | 119.4544 | 1.17 | 119.4544 | 0.000 % | 5.32 |
| SI-D1-C5-T3 | 119.3063 | 1.16 | 119.3063 | 0.000 % | 5.28 |
| SI-D2-C5-T3 | 119.4544 | 1.19 | 119.4544 | 0.000 % | 5.64 |
| SI-D1-C5-T4 | 119.3063 | 1.18 | 119.3063 | 0.000 % | 7.12 |
| SI-D2-C5-T4 | 119.4544 | 1.26 | 119.4544 | 0.000 % | 6.44 |

| Instance | $Z_1$ | Time | $Z_1$ | Gap | Time |
|---|---|---|---|---|---|
| SI-D1-C10-T2 | 119.9658 | 2.42 | 119.9658 | 0.000 % | 11.43 |
| SI-D2-C10-T2 | 120.1536 | 1.97 | 120.1536 | 0.000 % | 14.02 |
| SI-D1-C10-T3 | 119.9658 | 2.64 | 119.9658 | 0.000 % | 11.69 |
| SI-D2-C10-T3 | 120.1536 | 2.73 | 120.1536 | 0.000 % | 13.23 |
| SI-D1-C10-T4 | 119.9658 | 2.81 | 119.9658 | 0.000 % | 12.01 |
| SI-D2-C10-T4 | 120.1536 | 2.75 | 120.1536 | 0.000 % | 15.87 |
| SI-D1-C15-T2 | 286.1563 | 43.16 | 286.1563 | 0.000 % | 31.34 |
| SI-D2-C15-T2 | 286.1212 | 52.84 | 286.2178 | 0.034 % | 28.82 |
| SI-D1-C15-T3 | 231.3617 | 54.25 | 231.3617 | 0.000 % | 32.53 |
| SI-D2-C15-T3 | 230.9398 | 580.52 | 230.9398 | 0.000 % | 35.72 |
| SI-D1-C15-T4 | 231.3617 | 675.63 | 231.3617 | 0.000 % | 34.03 |
| SI-D2-C15-T4 | 230.9398 | 652.23 | 230.9398 | 0.000 % | 37.45 |
| SI-D1-C20-T2 | 295.9807 | 7200 | 295.9807 | 0.000 % | 88.74 |
| SI-D2-C20-T2 | 295.7998 | 7200 | 295.7854 | -0.005 % | 84.67 |
| SI-D1-C20-T3 | 296.8966 | 7200 | 296.7154 | -0.061 % | 95.62 |
| SI-D2-C20-T3 | 296.9625 | 7200 | 296.8275 | -0.045 % | 105.35 |
| SI-D1-C20-T4 | 296.7437 | 7200 | 296.7154 | -0.010 % | 82.38 |
| SI-D2-C20-T4 | 296.8345 | 7200 | 296.8275 | -0.002 % | 104.95 |
| SI-D1-C25-T2 | 298.0373 | 7200 | 297.5612 | -0.160 % | 137.42 |
| SI-D2-C25-T2 | 299.1454 | 7200 | 297.5629 | -0.434 % | 121.06 |
| SI-D1-C25-T3 | 299.3157 | 7200 | 298.2126 | -0.369 % | 132.14 |
| SI-D2-C25-T3 | 298.8595 | 7200 | 298.0089 | -0.285 % | 130.01 |
| SI-D1-C25-T4 | 298.3536 | 7200 | 298.2126 | -0.047 % | 141.77 |
| SI-D2-C25-T4 | 298.5121 | 7200 | 298.3093 | -0.068 % | 162.82 |
| Avg. | | 2949.36 | | -0.048 % | 56.67 |

As the results in Table 4 illustrate, our proposed solution algorithm establishes a good performance in solving small-size problems to optimality in a short time. For instance, sets of 5, 10, and 15 customers, the optimality of the results provided by CPLEX was guaranteed, and the computation time was shorter than two hours. Since the obtained gap is zero in these instances, the results provided by our proposed algorithm are also globally optimum. In instances with 20 and 25 customers, CPLEX was unable to reach the global optimal within two hours. In these instances, the gap of CPLEX to the linear relaxation of the objective function found in iterations of branch and bound was less than 0.7%. This implies that the obtained solutions by CPLEX were either globally optimal or very near to the global optimal. In these 12 instances, our solution algorithm achieved results equal to or smaller than those provided by CPLEX.

Table 5 presents the result of investigating the performance of our proposed algorithm that applies ALNS, hybridized with K-means and LS, together with B&P on benchmark instances of 50-200 customers. To this end, the objective function value of the Best Known Solution (BKS), which is found during the overall testing of the algorithm, is reported. Then, the performance of the proposed algorithm with and without its hybridization with K-means and LS is investigated, and related results, associated with the best solution of ten runs, including the objective function value ($Z_1$), its gap to BKS ($\Delta_{BKS}$), and run time ($t$), are reported.

Table 5. The results of experiments on newly generated medium and large-size benchmark instances

| Instance | BKS | ALNS+B&P(+ K-means +LS) | | | ALNS+B&P | | |
|---|---|---|---|---|---|---|---|
| | | $Z_1$ | $\Delta_{BKS}$ | $t$ (m) | $Z_1$ | $\Delta_{BKS}$ | $t$ (m) |
| MI-D3-C50-T5 | 395.9445 | 395.9445 | 0.000 % | 5.36 | 396.8985 | 0.241 % | 4.98 |
| MI-D4-C50-T5 | 394.9225 | 394.9225 | 0.000 % | 5.88 | 396.3533 | 0.362 % | 5.35 |
| MI-D3-C50-T10 | 369.9841 | 369.9841 | 0.000 % | 5.64 | 371.8903 | 0.515 % | 5.19 |
| MI-D4-C50-T10 | 368.2105 | 368.4828 | 0.074 % | 6.12 | 370.1435 | 0.525 % | 5.81 |
| MI-D3-C75-T5 | 437.5509 | 437.5509 | 0.000 % | 8.18 | 446.5644 | 2.056 % | 6.47 |
| MI-D4-C75-T5 | 435.0217 | 435.3131 | 0.067 % | 9.27 | 447.6808 | 2.913 % | 7.21 |
| MI-D3-C75-T10 | 398.1412 | 398.4636 | 0.081 % | 8.71 | 410.2845 | 3.048 % | 6.95 |
| MI-D4-C75-T10 | 396.8405 | 397.2135 | 0.094 % | 9.43 | 409.6187 | 3.215 % | 7.89 |
| LI-D3-C100-T5 | 498.3201 | 498.8184 | 0.102 % | 11.32 | 514.4656 | 3.235 % | 9.62 |
| LI-D4-C100-T5 | 492.0076 | 492.5981 | 0.117 % | 11.69 | 510.7531 | 3.812 % | 10.05 |
| LI-D3-C100-T10 | 471.0568 | 471.5749 | 0.105 % | 11.61 | 489.9461 | 4.007 % | 9.75 |
| LI-D4-C100-T10 | 468.1005 | 468.6622 | 0.122 % | 12.05 | 487.3862 | 4.116 % | 10.12 |
| LI-D3-C150-T5 | 614.2237 | 615.4275 | 0.196 % | 15.21 | 641.9374 | 4.512 % | 12.01 |
| LI-D4-C150-T5 | 601.5589 | 602.7561 | 0.199 % | 15.33 | 628.2801 | 4.442 % | 12.55 |
| LI-D3-C150-T10 | 583.0153 | 584.1171 | 0.189 % | 15.12 | 610.5394 | 4.721 % | 12.84 |
| LI-D4-C150-T10 | 580.1116 | 581.2892 | 0.203 % | 15.56 | 608.5892 | 4.909 % | 13.09 |
| LI-D3-C200-T5 | 807.4012 | 809.5891 | 0.271 % | 18.42 | 858.4045 | 6.317 % | 16.35 |
| LI-D4-C200-T5 | 804.3613 | 806.5006 | 0.266 % | 18.87 | 856.8133 | 6.521 % | 16.58 |
| LI-D3-C200-T10 | 785.0109 | 787.0264 | 0.257 % | 18.59 | 835.5792 | 6.442 % | 16.49 |
| LI-D4-C200-T10 | 783.0784 | 785.0591 | 0.253 % | 19.01 | 833.5943 | 6.451 % | 16.73 |
| Avg. | | | 0.131 % | 12.07 | | 3.618 % | 10.31 |

The observed results indicate that the intensified algorithm, which incorporates local search and K-means clustering in addition to adaptive large neighborhood search and branch and price, consistently outperforms the alternative algorithm in terms of solution quality. The average gap with the best-known solution (BKS) is approximately 3.5% lower for the intensified algorithm. Although the intensified algorithm requires an average execution time that is 1.5 minutes longer, the inclusion of K-means clustering helps mitigate this increase, resulting in a reasonably close execution time to the alternative algorithm. Furthermore, as the problem instances increase in size, the intensified algorithm exhibits a more significant performance advantage. This suggests that the additional time taken by the local search component is worthwhile, as it contributes to improved solutions for larger problem instances.

### 5.4. Case Study

Amsterdam's city center is a well-known hub, teeming with many hotels, restaurants, and cafés tightly clustered together. More than 1500 HoReCa spots exist in this area, requiring a huge distribution chain for their input food supplies. This section provides an optimal design for the distribution chain of restaurants and cafés located in Amsterdam's historical center.

The location of restaurants and cafés is derived from the municipality of Amsterdam's provided maps[1]. The municipality has also specified a set of potential locations for transshipment points and has classified them into poor, moderate, and spacious points based on the available space on the quay. Figure 2 illustrates the location of these TPs, together with restaurants and cafés.

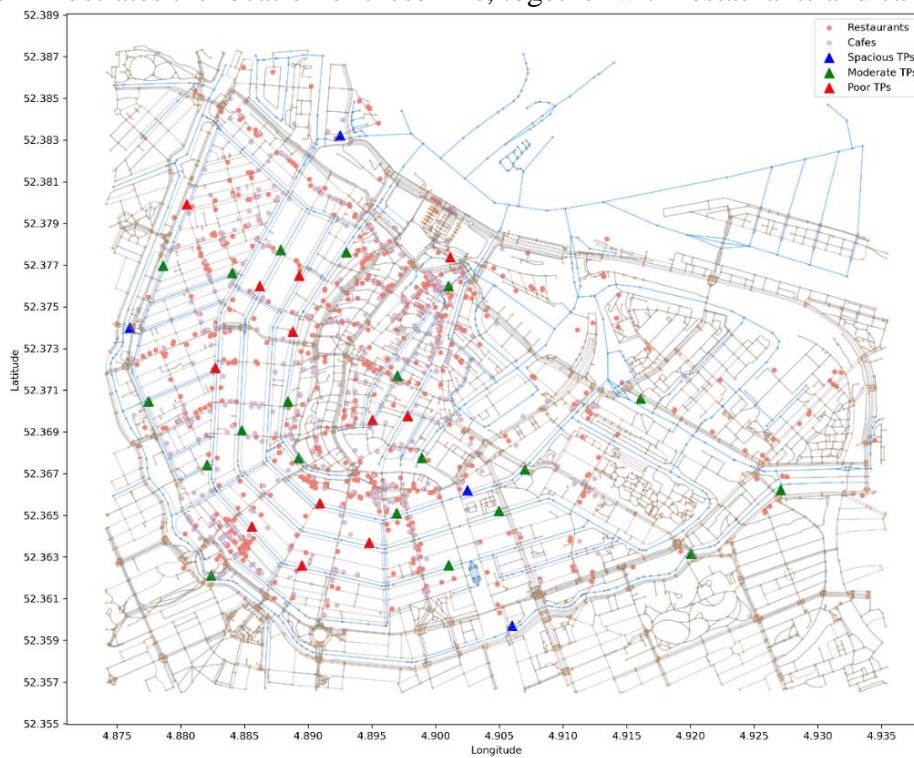

Figure 2. Restaurants and Cafés in the City Center of Amsterdam

---

[1] Available at https://maps.amsterdam.nl/functiekaart/

In order to estimate the daily demand of these businesses (in $m^3$), we have incorporated a Deep Neural Network (DNN) that works based on the input labels of business type (restaurant or café), weekday or weekend, season, longitude, and latitude. In order to obtain train data for our DNN, demand data were collected through field trip. More precisely, the city center was clustered into five location groups based on the longitude and latitude of the existing restaurants and cafés. Within each cluster, a number of restaurants and cafés was selected as visiting candidates. The number of visits was proportional to the number of all existing restaurants and cafés (based on intensity varying between 12%-20%). The daily demand during weekdays and weekend, and regarding four seasons, were estimated based on the collected information at visits.

*Results*

Figure 3 illustrates the result of solving the problem for the case of Amsterdam.

The results of our case study, as depicted in Figure 3, illustrates the establishment of 10 transshipment locations (TPs) strategically positioned within the city center, with a focus on densely populated areas. Among these TPs, the majority are categorized as moderate, while one is spacious, and two are considered poor in terms of size and capacity. A closer examination of the map reveals that the potential spacious TPs are not located in densely populated segments. This observation justifies their absence in the established set. The poor established TPs, on the other hand, are found only in segments where there is a lack of spacious or moderate TP candidates. To efficiently cater to the demands of the restaurants and cafés, a fleet of seven small and two medium vessels is employed, without the inclusion of any large vessels. Despite the potential cost-saving benefits associated with economies of scale, the use of large vessels was not feasible concerning the canal classes due to the limited width or depth of the canals within the city center. In the second echelon of our proposed system, 74 LEVs and 133 moving jacks are utilized to deliver the items to their final destinations. This indicates the appropriate location of established TPs has effectively reduced the number of LEVs required, thanks to the inclusion of moving jacks for serving points that are in close proximity to waterways.

A fundamental question in evaluating the efficiency of this waterway-based chain is if it can improve the distribution of the HoReCa demands. In order to answer this question, we compare the designed distribution chain with the one currently implemented in Amsterdam, for which trucks with a weight limit of 3500 kg deliver the items to demanded spots. This transforms the problem into a routing problem with time windows. Figure 4 compares these two distribution chains in terms of the total cost, the total number of applied road vehicles, their associated weight, and the average distance driven by each vehicle within the city center.

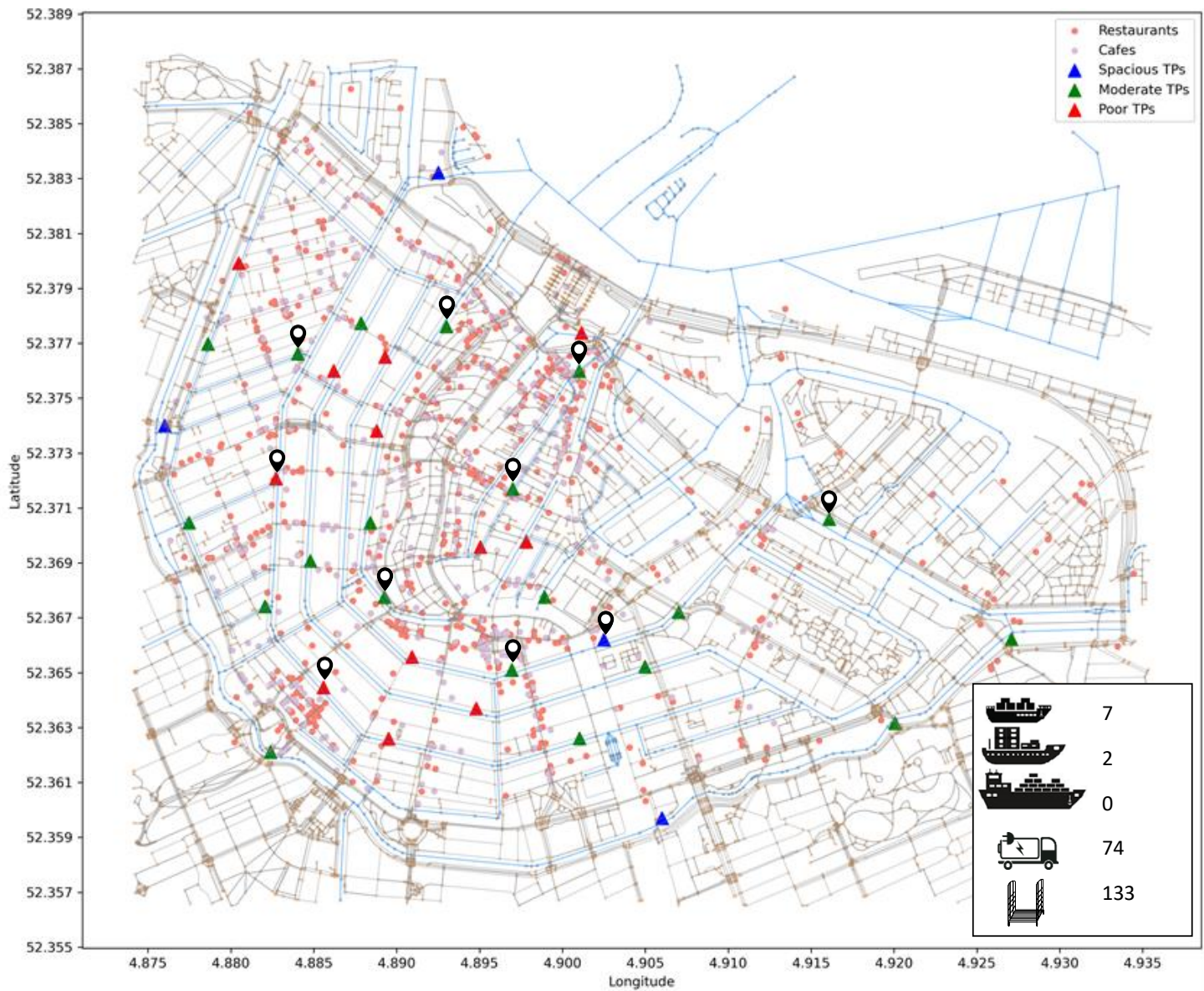

Figure 3. Case study Results

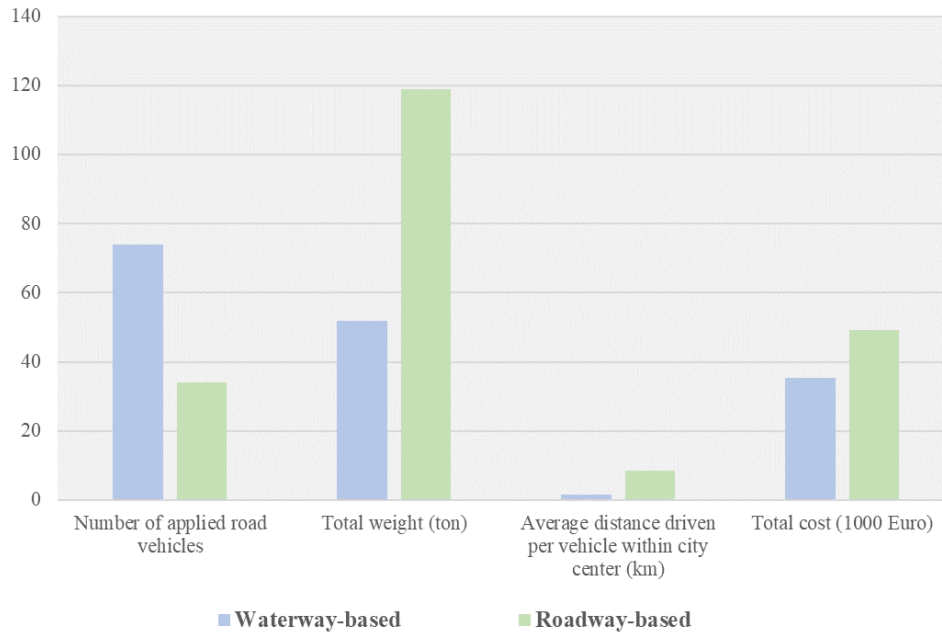

Figure 4. Comparison between waterway-based and roadway-based logistics chain

The analysis reveals that the waterway-based food distribution chain presents a noteworthy advantage in terms of total cost, potentially leading to cost savings of approximately 28% compared to the truck-based system. Despite the higher number of vehicles employed in the bi-modal setting, it is important to note that these vehicles are light vehicles, resulting in a 43% reduction in the total weight of vehicles driven within the city center. This not only offers the potential to preserve the lifetime of physical infrastructure such as quay walls and bridges but also indicates a more distributed and flexible delivery system.

Moreover, the waterway-based food chain demonstrates a significant 80% reduction in the average distance driven within the city center to serve the HoReCa spots compared to the truck-based system. This reduction in distance traveled has the potential to alleviate traffic congestion, improve efficiency in terms of time and fuel consumption, and contribute to decreased emissions. Furthermore, since the bi-modal setting utilizes electric vehicles, it offers a substantial reduction in carbon emissions as well.

In conclusion, the results strongly suggest that implementing the waterway-based distribution chain has the potential to enhance the efficiency of HoReCa demands in Amsterdam greatly. The advantages encompass lower total cost, a more distributed fleet of lighter vehicles, a significant reduction in average distance driven within the city center, and a notable decrease in exhaust emissions.

### 5.5. Sensitivity Analysis

In this sub-section, we will investigate the impact of different input parameters on our designed distribution chain. In order to be able to track this impact on all decision variables, a medium instance with three depots, five transshipment locations, and 50 customers is selected.

The parameters associated with the period equivalent establishment cost and travel cost are the two cost factors characterizing the economic competency of our designed waterway-based distribution chain. In order to investigate their impact on total cost, number of established TPs, number of applied vessels, and number of applied LEVs, sensitivity analyses are carried out on these parameters, and the results are provided in Table 6.

Table 6. The impact of changes on different cost parameters

| Parameter | Changes (%) | Total Cost | $n_{TP}$ | $n_V$ | $n_{LEV}$ |
|---|---|---|---|---|---|
| $FC_i$ | -75 % | 129.85 | 5 | 3 | 21 |
|  | -50 % | 230.55 | 5 | 3 | 21 |
|  | -25 % | 308.07 | 4 | 3 | 26 |
|  | 0 | 395.94 | 4 | 3 | 26 |
|  | +25 % | 486.24 | 4 | 3 | 26 |
|  | +50 % | 575.24 | 3 | 3 | 32 |
|  | +75 % | 664.24 | 3 | 3 | 32 |
| $C_{ijk}^{I}$ | -75 % | 362.54 | 4 | 4 | 27 |
|  | -50 % | 379.88 | 4 | 3 | 26 |
|  | -25 % | 391.51 | 4 | 3 | 26 |
|  | 0 | 395.94 | 4 | 3 | 26 |
|  | +25 % | 405.51 | 4 | 3 | 26 |
|  | +50 % | 414.54 | 4 | 3 | 26 |
|  | +75 % | 442.38 | 4 | 3 | 26 |
| $C_{ijk}^{II}$ | -75 % | 390.99 | 4 | 3 | 28 |
|  | -50 % | 392.64 | 4 | 3 | 27 |
|  | -25 % | 394.29 | 4 | 3 | 27 |
|  | 0 | 395.94 | 4 | 3 | 26 |
|  | +25 % | 397.59 | 4 | 3 | 26 |
|  | +50 % | 399.23 | 4 | 3 | 26 |
|  | +75 % | 400.89 | 4 | 3 | 25 |

The results indicate that changes in transshipment location establishment costs have a significantly larger influence on total costs compared to variations in first and second echelon transportation costs. Specifically, by reducing the establishment costs, the number of established locations may increase, leading to a decreased reliance on LEVs for transportation. This reduction in LEV usage was attributed to the improved efficiency achieved through the utilization of moving jacks for item delivery. This indicates that by effectively lowering this cost, one not only benefits from reductions in cost but also can decrease the road traffic load by fewer applied vehicles. Conversely, the variations in first and, specifically, second echelon transportation costs were found to have a relatively minor impact on overall costs. These findings underscore the importance of effectively managing and optimizing transshipment location establishment costs as a key strategy for achieving cost efficiencies in transportation operations, while lower attention can be devoted to the second echelon transportation.

**5.6. Discussion**

In our study, we conducted a thorough numerical analysis to explore the performance and potential benefits of a proposed algorithm and the implementation of a waterway-based distribution chain in Amsterdam. The results provide an intriguing story of improved efficiency and cost savings.

Our findings revealed that our algorithm consistently outperformed the previous approaches, indicating its effectiveness in solving the existing benchmarks and further our developed benchmark instance. This suggests that investing in algorithm development and optimization can lead to improved performance and efficiency in business operations. Motivated by these promising results, we delved deeper into the advantages of implementing a waterway-based distribution chain. Through a comprehensive case study, we compared this system to a traditional truck-based approach. The analysis considered various factors such as total cost, vehicle weight, distance traveled within the city center and environmental impact.

The outcomes of our case study form a compelling narrative. Implementing the waterway-based distribution chain yielded significant cost savings of approximately 28% compared to the truck-based system. This financial advantage was accompanied by a more distributed fleet of lighter vehicles. By utilizing lighter vehicles, we reduced the total weight of vehicles driven within the city center by 43%. This not only benefits the infrastructure's longevity but also allows for a more flexible and adaptable delivery system that can be of particular interest to city authorities and policymakers. Moreover, the waterway-based food chain showcased an impressive 80% reduction in the average distance traveled within the city center, offering the potential to alleviate traffic congestion, enhance efficiency in terms of time and fuel consumption, and contribute to decreased emissions. In addition, the adoption of electric vehicles in the bi-modal setting further reduced carbon emissions, reinforcing the environmental advantages of the waterway-based distribution chain.

Our sensitivity analysis shows that effectively managing and optimizing transshipment location establishment costs is crucial for achieving cost efficiencies in transportation operations. Reducing these costs not only leads to overall cost savings but also decreases reliance on LEVs and reduces road traffic congestion. While first and second echelon transportation costs have a relatively minor impact, organizations should prioritize optimizing transshipment location establishment costs as a key strategy for cost efficiency.

In conclusion, our numerical analysis presents a compelling case for both the proposed algorithm and the implementation of a waterway-based distribution chain in Amsterdam. The algorithm showcased superior performance, while the waterway-based system offered notable cost savings, reduced traffic congestion, improved efficiency, and decreased emissions. These insights provide valuable guidance for managers seeking to enhance operational efficiency, reduce costs, and contribute to sustainable transportation practices.

## 6. Conclusion

This study addresses the need for efficient urban logistics in Amsterdam by proposing an optimal urban logistic network that integrates urban waterways and last-mile delivery. The research highlights the untapped potential of inland waterways and the benefits they can offer in terms of capacity and addressing logistical challenges in the city center. In this respect, we propose an optimal urban logistic network for Amsterdam, integrating urban waterways and last-mile delivery via road transportation. The problem is formulated as a two-echelon location routing problem with time windows, and a hybrid solution approach is developed to solve it effectively.

The proposed algorithm consistently outperforms existing ones, demonstrating its effectiveness in solving existing benchmarks and newly developed instances. Through a comprehensive case study, the advantages of implementing a waterway-based distribution chain were assessed. The waterway-based chain showcased significant cost savings of approximately 28% compared to the traditional truck-based system. Additionally, the adoption of lighter vehicles led to a 43% reduction in total vehicle weight within the city center, enhancing infrastructure longevity and enabling a more flexible delivery system.

Furthermore, the waterway-based chain demonstrated an impressive 80% reduction in the average distance traveled within the city center, which has the potential to alleviate traffic congestion, improve efficiency in terms of time and fuel consumption, and contribute to decreased emissions. The incorporation of electric vehicles further reduced carbon emissions, highlighting the environmental advantages of the proposed system. The sensitivity analysis emphasized the importance of effectively managing and optimizing transshipment location establishment costs as a key strategy for achieving cost efficiencies and reducing reliance on LEVs and road traffic congestion. By drawing insights from Amsterdam's experience and embracing innovative approaches, cities around the world can endeavor to discover sustainable solutions for their urban logistics challenges.

Overall, this study provides valuable insights and guidance for managers in their pursuit of enhancing operational efficiency, reducing costs, and promoting sustainable transportation practices. However, further analysis is necessary to fully assess the viability and potential benefits of implementing the waterway-based chain, including considerations of infrastructural limitations and canal features. It is also important to study the impact of such a modal shift on the traffic over water, and increase in the propeller wash that in turn may lead to further deterioration of bed level and quay walls. A digital twin of Amsterdam's city center canals can be developed in that respect, which would illustrate the impact of changes in network design. Establishing a feedback loop between the optimization and the digital twin is another interesting way forward. The digital twin can get insights into the quality of solutions obtained by optimization and provide it with guidance on specific promising or forbidden solution space directions. Via such a setting, one can benefit both from design capability of the optimization and what-if scenario analysis of the digital twin.